\documentclass[12pt]{article}
\usepackage[T1]{fontenc}
\usepackage{lmodern}
\usepackage{amsmath,amssymb,amsthm,mathtools}
\usepackage{microtype}
\usepackage[top=27mm, bottom=27mm, left=22mm, right=22mm]{geometry}
\usepackage{xcolor}
\usepackage{tikz}
\usepackage{placeins}
\usepackage{cite}
\usepackage[colorlinks=true,linkcolor=blue!50!black,citecolor=blue!50!black,urlcolor=blue!50!black]{hyperref}

\newtheorem{theorem}{Theorem}[section]
\newtheorem{lemma}[theorem]{Lemma}
\newtheorem{proposition}[theorem]{Proposition}
\newtheorem{corollary}[theorem]{Corollary}
\theoremstyle{remark}
\newtheorem{remark}[theorem]{Remark}
\newcommand{\F}{\mathbb{F}}

\title{A near-quadratic lower bound for the Brown--Erd\H{o}s--S\'os $(9,6)$-problem}
\author{Xiangliang~Kong\thanks{Xiangliang~Kong is with the State Key Laboratory of Mathematical Sciences, Academy of Mathematics and Systems Science, Chinese Academy of Sciences, Beijing 100190, China. Email: \mbox{\texttt{xlkong@amss.ac.cn}}}
\and Chong~Shangguan\thanks{Chong~Shangguan is with the Research Center for Mathematics and Interdisciplinary Sciences, Shandong University, Qingdao 266237, China, and the Frontiers Science Center for Nonlinear Expectations, Ministry of Education, Qingdao 266237, China. Email: \mbox{\texttt{theoreming@163.com}}}
\and Xin~Wang\thanks{Xin~Wang is with the Department of Mathematics, Soochow University, Suzhou 215005, Jiangsu, China. Email: \mbox{\texttt{xinw@suda.edu.cn}}}}
\date{}

\begin{document}
\maketitle

\begin{abstract}
For integers $e\ge3$, let $f_3(n,e+3,e)$ denote the maximum number of edges in an $n$-vertex $3$-uniform hypergraph containing no $e$ distinct edges spanning at most $e+3$ vertices. The lower-bound part of a longstanding conjecture of Erd\H{o}s, Frankl, and R\"odl from 1986 asserts that $f_3(n,e+3,e)\ge n^{2-o(1)}$ for every fixed $e\ge3$. We prove that $f_3(n,9,6)\ge n^2\exp(-O((\log n)^{2/3}))$, resolving the first open case not covered by the classical construction of Ruzsa and Szemer\'edi in 1976, which established $f_3(n,6,3)\ge n^{2-o(1)}$. We also establish near-quadratic lower bounds for $f_3(n,e+3,e)$ with $e=9,10,11$. We further extend the result for $e=6$ to every fixed uniformity $r\ge4$.
\end{abstract}

\section{Introduction}

For integers $n,e\ge3$, let $f_3(n,e+3,e)$ denote the maximum number of edges in an $n$-vertex $3$-uniform hypergraph containing no $e$ distinct edges spanning at most $e+3$ vertices. The Brown--Erd\H{o}s--S\'os conjecture from 1973 \cite{BrownErdosSos1973} asserts that $f_3(n,e+3,e)=o(n^2)$ as $n\to\infty$ for every fixed $e\ge3$. The only settled case is $e=3$, established by the celebrated $(6,3)$-theorem of Ruzsa and Szemer\'edi \cite{RuzsaSzemeredi1978}. Their construction, based on Behrend's progression-free sets \cite{Behrend1946}, also shows that $f_3(n,6,3)>n^{2-o(1)}$ (see also \cite{GowersJanzer2021} for a different geometric construction). Despite much effort \cite{SarkozySelkow2005,SolymosiSolymosi2017,ConlonGishbolinerLevanzovShapira2023,JanzerMethukuMilojevicSudakov2025,ShapiraTyomkyn2025}, the conjecture remains open for every $e\ge4$.

The complementary near-quadratic lower-bound conjecture also has a long history. More generally, for integers $r,e\ge3$, an $r$-uniform hypergraph is \emph{$(e(r-2)+3,e)$-free} if no $e$ distinct edges span at most $e(r-2)+3$ vertices. Let $f_r(n,e(r-2)+3,e)$ denote the maximum number of edges in an $n$-vertex $(e(r-2)+3,e)$-free $r$-uniform hypergraph. In 1986, Erd\H{o}s, Frankl, and R\"odl \cite[Problem~6.2]{ErdosFranklRodl1986} asked whether $n^{2-o(1)}\le f_r(n,e(r-2)+3,e)=o(n^2)$ holds for every fixed $r,e\ge3$ as $n\rightarrow\infty$. Alon and Shapira \cite{AlonShapira2006} subsequently formulated a further generalization; see also \cite{ConlonGishbolinerLevanzovShapira2023,GeShangguan2021,ShapiraICM2022}. The conjectured lower bound would show that the quadratic exponent in the Brown--Erd\H{o}s--S\'os conjecture is best possible.

The probabilistic construction of Brown, Erd\H{o}s, and S\'os \cite{BrownErdosSos1973r} gives $f_r(n,e(r-2)+3,e)=\Omega(n^{2-1/(e-1)})$ for every fixed $r,e\ge3$. Shangguan and Tamo \cite{ShangguanTamo2020} improved this to $\Omega(n^{2-1/(e-1)}(\log n)^{1/(e-1)})$. The conjectured near-quadratic lower bound for $e=3$ and every fixed $r\ge3$ was established by Erd\H{o}s, Frankl, and R\"odl \cite{ErdosFranklRodl1986}, extending the result of Ruzsa and Szemer\'edi \cite{RuzsaSzemeredi1978}. Ge and Shangguan \cite{GeShangguan2021} subsequently established it for $e=4,5,7,8$ and every fixed $r\ge3$. For $r=3$, the cases $e=4,5$ follow directly from the Ruzsa--Szemer\'edi $(6,3)$ construction \cite{AlonShapira2006}.

The first nontrivial case left open by these results was therefore $e=6$, even for $r=3$. In particular, the question whether $f_3(n,9,6)\ge n^{2-o(1)}$ was explicitly raised in several works \cite{GeShangguan2021,GishbolinerShapira2022,ConlonGishbolinerLevanzovShapira2023} and discussed in Shapira's ICM survey \cite[Section~2.2]{ShapiraICM2022}. Specializing the result of Shangguan and Tamo \cite{ShangguanTamo2020} gives $f_3(n,9,6)=\Omega(n^{9/5}(\log n)^{1/5})$, which falls short of the conjectured near-quadratic bound.

A $3$-uniform hypergraph is \emph{linear} if any two distinct edges share at most one vertex. A \emph{loose triangle} consists of three edges such that each pair intersects in exactly one vertex and the three intersection vertices are distinct. A \emph{$3\times3$ grid} consists of three pairwise disjoint row edges and three pairwise disjoint column edges such that every row edge intersects every column edge. As observed in \cite{GeShangguan2021,GishbolinerShapira2022}, every linear $3$-uniform hypergraph with six edges on at most nine vertices either contains a loose triangle or is a $3\times3$ grid. Consequently, any linear $3$-uniform hypergraph containing neither configuration is $(9,6)$-free.

F\"uredi and Ruszink\'o \cite{FurediRuszinko2013} constructed linear $3$-uniform hypergraphs containing neither loose triangles nor grids with $n^{8/5-o(1)}$ edges, and Ge and Shangguan \cite{GeShangguan2021} improved this bound to $\Omega(n^{8/5})$. In contrast, excluding grids alone allows quadratic density. Gishboliner and Shapira \cite{GishbolinerShapira2022} constructed linear grid-free $3$-uniform hypergraphs with $(1/16-o(1))n^2$ edges, and Solymosi \cite{Solymosi2026} improved this to $(1/12-o(1))n^2$; see also \cite{Pohoata2026} for a different construction with a weaker bound when $r=3$. These grid-free constructions do not guarantee the absence of loose triangles and therefore do not establish $(9,6)$-freeness. Thus the challenge is to exclude both configurations while retaining $n^{2-o(1)}$ edges.

\subsection{Main results}

Our first result gives a linear tripartite $3$-uniform hypergraph with $n^{2-o(1)}$ edges that is $(9,6)$-free.

\begin{theorem}\label{thm:main}
There is an absolute constant $C>0$ such that, for every sufficiently large $n$, there exists a linear tripartite $3$-uniform hypergraph on $n$ vertices that has at least $n^2\exp(-C(\log n)^{2/3})$ edges, contains neither a loose triangle nor a $3\times3$ grid, and is $(9,6)$-free. In particular,
\begin{equation*}
f_3(n,9,6)\ge n^2\exp\bigl(-C(\log n)^{2/3}\bigr).
\end{equation*}
\end{theorem}

\begin{remark}\label{rem:further}
The same construction is in fact $(e+3,e)$-free for every integer $3\le e\le11$. Consequently, the same lower bound holds for $f_3(n,e+3,e)$ with $e=9,10,11$. The additional properties are established in Appendix~\ref{app:further}.
\end{remark}

The conjectured lower bound $f_3(n,e+3,e)\ge n^{2-o(1)}$ was previously known for $e=3,4,5,7,8$ \cite{RuzsaSzemeredi1978,AlonShapira2006,GeShangguan2021}. Theorem~\ref{thm:main} establishes it for $e=6$, while Remark~\ref{rem:further} gives the remaining cases $e=9,10,11$. 

We next extend the result for $e=6$ to every fixed $r\ge4$.

\begin{theorem}\label{thm:higher}
For every fixed integer $r\ge4$, there is a constant $C_r>0$ such that, for every sufficiently large $n$, there exists a linear $r$-partite $r$-uniform $(6r-9,6)$-free hypergraph on $n$ vertices with at least
\begin{equation*}
n^2\exp\bigl(-C_r(\log n)^{\theta_r}\bigr)
\quad\text{edges},\qquad \theta_r=\frac{2^{r-1}}{2^{r-1}+1}.
\end{equation*}
Consequently, $f_r(n,6r-9,6)\ge n^{2-o(1)}$ for every fixed $r\ge3$.
\end{theorem}

We also obtain an application to extremal set theory. A hypergraph is \emph{$3$-union-free} if any two distinct subfamilies of at most three edges have different unions. Let $U_3(n,3)$ denote the maximum number of edges in an $n$-vertex $3$-uniform hypergraph that is $3$-union-free. F\"uredi and Ruszink\'o~\cite{FurediRuszinko2013} proved that $U_3(n,3)=\Omega(n^{5/3})$. More recently, Liu, Shangguan, and Zhang~\cite{LiuShangguanZhang2026} obtained sharp asymptotic bounds for many uniform union-free hypergraphs, but left open whether $U_3(n,3)\ge n^{2-o(1)}$.

As shown by F\"uredi and Ruszink\'o~\cite[Proof of Corollary~1.9, Eq.~(12)]{FurediRuszinko2013}, every linear $3$-uniform hypergraph containing neither a loose triangle nor a $3\times3$ grid is $3$-union-free. Thus Theorem~\ref{thm:main} yields the following corollary.

\begin{corollary}
There is an absolute constant $C>0$ such that, for every sufficiently large $n$, there exists a linear $3$-uniform $3$-union-free hypergraph on $n$ vertices with at least $n^2\exp(-C(\log n)^{2/3})$ edges. In particular, $U_3(n,3)>n^{2-o(1)}$.
\end{corollary}

\subsection{Proof overview}

Earlier constructions of F\"uredi and Ruszink\'o~\cite{FurediRuszinko2013} and Ge and Shangguan~\cite{GeShangguan2021} considered $3$-partite $3$-uniform hypergraphs with edges of the form $(y,y+m,y+2m)$ and $(y,y+m,y+3m)$, respectively, where $y\in G$ and $m\in S\subseteq G$ for a suitable finite abelian group $G$. Additive restrictions on $S$ are used to exclude loose triangles and grids, but these constructions did not achieve $n^{2-o(1)}$ edges. Indeed, one possible approach to achieving such a bound would be to find a set $S\subseteq G$ of size $|G|^{1-o(1)}$ such that the equation $2\alpha+2\beta=3\gamma+\delta$ has only the trivial solutions $\alpha=\beta=\gamma=\delta$ in $S$. This is a difficult problem studied by Ruzsa~\cite{Ruzsa1993}. In fact, the existence of such sets in the integers would give $3$-uniform hypergraphs with $n^{2-o(1)}$ edges in which no five edges span at most nine vertices, contradicting a conjecture of Gowers and Long~\cite{GowersLong2021}. Recently, Hametner and Solymosi~\cite{HametnerSolymosi2026} obtained the best known lower bound $n^{\alpha-o(1)}$ for the size of such a subset of $\{1,\ldots,n\}$, where $\alpha=\log 15/(3\log 5)>0.5608687$.

By contrast, Gishboliner and Shapira~\cite{GishbolinerShapira2022} considered triples of the form $(x,x+a,xa)$, where $x$ and $a$ range over suitable subsets of a field. The quadratic term $xa$ leads to an algebraic identity that rules out grids, allowing them to construct linear grid-free $3$-graphs with $\Omega(n^2)$ edges. Their construction, however, does not in general exclude loose triangles.

Our construction represents vertices by vectors in $\F_p^d$, where $p$ is prime and $d\ge2$, and uses a cubic term involving consecutive coordinates to combine the progression-free and algebraic arguments needed to exclude both configurations. More precisely, let $S\subseteq\F_p$ be a set containing no nontrivial three-term arithmetic progression. We take three disjoint copies of $\F_p^d$ as vertex parts and define edges $(x,y,z)$ by $z_1=x_1+y_1$ and $z_j=x_j+y_j+x_{j-1}y_{j-1}^{\,2}$ for $2\le j\le d$, where $y_j-x_j\in S$ for $1\le j<d$ and $y_d=x_d$. The key idea is to delay the cubic term by one coordinate. When analyzing a loose triangle or a grid, we consider the first coordinate at which the corresponding $x$-values or $y$-values differ. At this coordinate, the cubic terms agree, leaving additive relations. For a loose triangle, these relations yield a forbidden three-term arithmetic progression in $S$ (see Lemma~\ref{lem:triangle}). For a grid, they force three coordinate values to be pairwise distinct, while the cubic terms at the next coordinate give a vanishing Vandermonde determinant, a contradiction (see Lemma~\ref{lem:grid}).

To extend the result for $e=6$ to higher uniformities, we use a related construction in which the vertex parts are indexed by distinct field elements $b$. In the vertex part indexed by $b$, the $j$th coordinate contains the linear term $u_j+bm_j$, with an additional cubic term $b^2m_{j-1}^{\,3}$ for $j\ge2$ (see \eqref{higher:eq:map}). The details are given in Section~\ref{sec:higher}.

The rest of this paper is organized as follows. Section~\ref{sec:construction} presents the three-uniform construction and proves Theorem~\ref{thm:main}. Section~\ref{sec:higher} extends the construction to higher uniformities and proves Theorem~\ref{thm:higher}. Appendix~\ref{app:further} establishes the additional properties stated in Remark~\ref{rem:further}. All logarithms are natural.

\section{The three-uniform case}\label{sec:construction}
\subsection{The construction and basic properties}
Behrend~\cite{Behrend1946} proved that every sufficiently large integer $N$ admits a subset of $\{0,\ldots,N-1\}$ of size at least $N\exp(-O(\sqrt{\log N}))$ containing no nontrivial three-term arithmetic progression, that is, no three distinct elements $a,b,c$ satisfying $a+c=2b$. Consequently, there is an absolute constant $C_0>0$ such that, for every sufficiently large prime $p$, there exists a set $S\subseteq\F_p$ satisfying
\begin{equation}\label{eq:behrend}
|S|\ge p\exp\bigl(-C_0\sqrt{\log p}\bigr),\qquad
a+c=2b\Longrightarrow a=b=c\quad(a,b,c\in S).
\end{equation}
Indeed, choose a Behrend set in $\{0,\ldots,\lfloor(p-1)/3\rfloor\}$ and regard it as a subset of $\F_p$. Since $|2b-a-c|\le 2\lfloor(p-1)/3\rfloor<p$ for elements of this interval, reduction modulo $p$ introduces no nontrivial three-term arithmetic progressions. The constant factor in the interval length is absorbed into $C_0$.

Fix a prime $p>3$, an integer $d\ge2$, and a set $S\subseteq\F_p$ containing no nontrivial three-term arithmetic progression. Let $V_1,V_2,V_3$ be disjoint copies of $\F_p^d$. For $x,y\in\F_p^d$, define $z=z(x,y)\in\F_p^d$ by
\begin{equation}\label{eq:map}
z_1=x_1+y_1,\qquad z_j=x_j+y_j+x_{j-1}y_{j-1}^{\,2}\quad(2\le j\le d).
\end{equation}
Let $H=H(p,d,S)$ be the tripartite $3$-uniform hypergraph with vertex parts $V_1,V_2,V_3$ and edges $(x,y,z(x,y))$ for all $x,y\in\F_p^d$ satisfying
\begin{equation}\label{eq:admissible}
y_j-x_j\in S\quad(1\le j<d),\qquad y_d=x_d.
\end{equation}

Below we establish some basic properties of $H$.

\begin{lemma}\label{lem:linear-count}
The hypergraph $H$ is linear, has $3p^d$ vertices, and has exactly $p^d|S|^{d-1}$ edges.
\end{lemma}

\begin{proof}
For each $x\in\F_p^d$, the conditions in \eqref{eq:admissible} allow exactly $|S|^{d-1}$ choices of $y$. Since distinct pairs $(x,y)$ determine distinct edges, we have $|E(H)|=p^d|S|^{d-1}$. Also, $|V(H)|=3p^d$.

To prove linearity, it suffices to show that any two vertices from different parts determine at most one edge. Given $x,y$, the vertex $z$ is uniquely determined by \eqref{eq:map}. Given $x,z$, we can recover $y$ successively from $y_1=z_1-x_1$ and $y_j=z_j-x_j-x_{j-1}y_{j-1}^2$ for $2\le j\le d$. Similarly, given $y,z$, we can recover $x$ successively from $x_1=z_1-y_1$ and $x_j=z_j-y_j-x_{j-1}y_{j-1}^2$ for $2\le j\le d$. Hence $H$ is linear.
\end{proof}

\subsection{Excluding triangles and grids}
The following classification appears in Gishboliner and Shapira \cite{GishbolinerShapira2022}; see also Ge and Shangguan \cite{GeShangguan2021} for a related multipartite result.

\begin{lemma}\label{lem:classification}
Let $F$ be a linear $3$-uniform hypergraph with six edges and at most nine vertices. Then $F$ contains a loose triangle or is a $3\times3$ grid.
\end{lemma}

\begin{proof}
Suppose that $F$ contains no loose triangle. If three edges share a vertex, their union $U$ has seven vertices. Any other edge meets $U$ in at most one vertex, since two intersections would violate linearity or form a loose triangle. As $F$ has at most nine vertices, each of the remaining three edges must contain the same two vertices outside $U$, contradicting linearity. Thus every vertex has degree at most two.

Since $F$ has six edges and at most nine vertices, it has exactly nine vertices, all of degree two. Let $G$ be the graph on the six edges of $F$, with adjacency given by intersection. By linearity, $G$ is $3$-regular; it is triangle-free because a triangle would give a loose triangle in $F$, as every vertex of $F$ has degree two. The only triangle-free $3$-regular graph on six vertices is $K_{3,3}$: the three neighbors of any vertex are independent and must be adjacent to both remaining vertices. Its two vertex classes give two triples of pairwise disjoint edges of $F$, with every opposite pair intersecting. Hence $F$ is a $3\times3$ grid.
\end{proof}

We first use the progression-free property of $S$ to exclude loose triangles.

\begin{lemma}\label{lem:triangle}
The hypergraph $H$ contains no loose triangle.
\end{lemma}

\begin{proof}
Suppose that $H$ contains a loose triangle. Since any two of its three intersection vertices lie in a common edge, they must belong to different vertex parts. Thus we may write the three edges as $(x,y,z)$, $(x,y',z')$, and $(x',y,z')$, where $x\ne x'$ and $y\ne y'$. In particular, $z(x,y')=z(x',y)$.

By \eqref{eq:admissible}, we have $x_d=y_d=y'_d=x'_d$. Let $t<d$ be the least coordinate for which $x_t\ne x'_t$ or $y_t\ne y'_t$. Since $x_j=x'_j$ and $y_j=y'_j$ for $j<t$, the cubic terms in coordinate $t$ of $z(x,y')$ and $z(x',y)$ agree (and are absent if $t=1$). Hence $x_t+y'_t=x'_t+y_t$, giving $\delta=x'_t-x_t=y'_t-y_t\ne0$ by the choice of $t$.

Put $a=y_t-x_t$. Since all three edges satisfy \eqref{eq:admissible}, we have $\{y_t-x_t,y'_t-x_t,y_t-x'_t\}=\{a,a+\delta,a-\delta\}\subseteq S$, contradicting the progression-free property of $S$.
\end{proof}

We next turn to grids. The following standard representation of a tripartite grid appears in earlier work \cite{FurediRuszinko2013,GishbolinerShapira2022}. We include a short proof via a proper edge-coloring of $K_{3,3}$, as in Ge and Shangguan \cite[Proof of Theorem~5.6]{GeShangguan2021}.

\begin{lemma}\label{lem:grid-normal}
Suppose that six edges of a linear tripartite $3$-uniform hypergraph with vertex parts $V_1,V_2,V_3$ form a $3\times3$ grid. Then there exist distinct vertices $X_1,X_2,X_3\in V_1$, distinct vertices $Y_1,Y_2,Y_3\in V_2$, and distinct vertices $Z_{12},Z_{13},Z_{23}\in V_3$ such that the six edges are
\begin{equation}\label{eq:grid-normal}
(X_i,Y_j,Z_{ij})\qquad(1\le i,j\le3,\ i\ne j),\qquad Z_{ji}=Z_{ij}.
\end{equation}
\end{lemma}

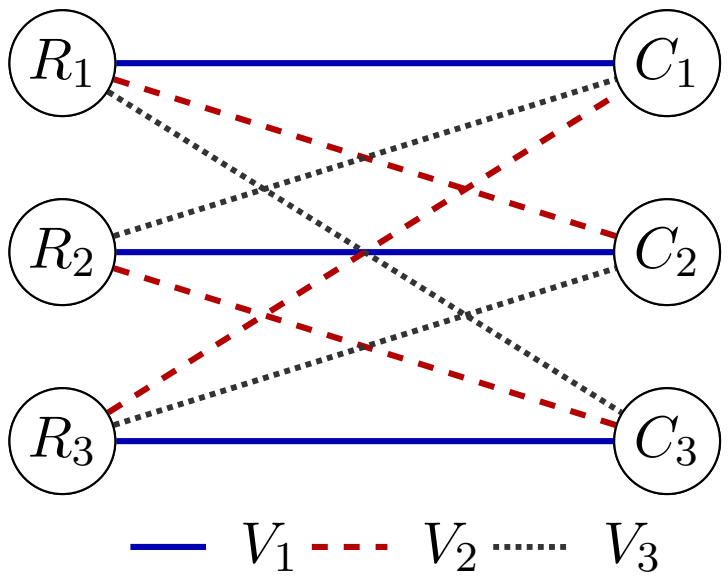
\begin{figure}[htbp]
\centering
\begin{minipage}[c]{0.40\textwidth}
\centering
$\displaystyle
\begin{array}{c|ccc}
 & C_1 & C_2 & C_3 \\ \hline
R_1 & X_1 & Y_3 & Z_{13} \\
R_2 & Z_{12} & X_2 & Y_1 \\
R_3 & Y_2 & Z_{23} & X_3
\end{array}$
\end{minipage}\hfill
\begin{minipage}[c]{0.58\textwidth}
\centering
\begin{tikzpicture}[x=1cm,y=1cm,every node/.style={font=\small}]
\tikzset{gridnode/.style={draw,circle,fill=white,minimum size=7mm,inner sep=1pt}}
\foreach \i/\y in {1/1.25,2/0,3/-1.25} {
  \node[gridnode] (R\i) at (0,\y) {$R_{\i}$};
  \node[gridnode] (C\i) at (4,\y) {$C_{\i}$};
}
\draw[blue!70!black,line width=1.15pt]
  (R1)--(C1) (R2)--(C2) (R3)--(C3);
\draw[red!70!black,dashed,line width=1.15pt]
  (R1)--(C2) (R2)--(C3) (R3)--(C1);
\draw[black!80,densely dotted,line width=1.1pt]
  (R1)--(C3) (R2)--(C1) (R3)--(C2);
\foreach \i/\y in {1/1.25,2/0,3/-1.25} {
  \node[gridnode] at (0,\y) {$R_{\i}$};
  \node[gridnode] at (4,\y) {$C_{\i}$};
}
\draw[blue!70!black,line width=1.15pt]
  (0.45,-1.95)--(0.95,-1.95);
\node[anchor=west] at (1.03,-1.95) {$V_1$};
\draw[red!70!black,dashed,line width=1.15pt]
  (1.65,-1.95)--(2.15,-1.95);
\node[anchor=west] at (2.23,-1.95) {$V_2$};
\draw[black!80,densely dotted,line width=1.1pt]
  (2.85,-1.95)--(3.35,-1.95);
\node[anchor=west] at (3.43,-1.95) {$V_3$};
\end{tikzpicture}
\end{minipage}
\caption{The labeling of the nine grid intersections (left) and the corresponding proper $3$-edge-coloring of $K_{3,3}$ (right). Each edge $R_iC_j$ in the graph represents the vertex $R_i\cap C_j$, and its color indicates the vertex part containing this intersection.}
\label{fig:grid-coloring}
\end{figure}

\begin{proof}
Let $R_1,R_2,R_3$ and $C_1,C_2,C_3$ be the row and column edges of the grid. Consider the complete bipartite graph $K_{3,3}$ with vertex classes $\{R_1,R_2,R_3\}$ and $\{C_1,C_2,C_3\}$. Color each edge $R_iC_j$ with $a\in\{1,2,3\}$ if the unique vertex in $R_i\cap C_j$ lies in $V_a$. Since each row and column edge contains exactly one vertex from each part, this is a proper $3$-edge-coloring of $K_{3,3}$, so each color class is a perfect matching.

After relabeling, we may assume that the edges of color $1$ are $R_iC_i$ for $1\le i\le3$. The edges of color $2$ form a perfect matching disjoint from these three edges. Since there are only two such matchings, equivalent under relabeling, we may assume that these edges are $R_1C_2$, $R_2C_3$, and $R_3C_1$. The remaining three edges then have color $3$. Thus the nine intersections can be labeled as in Figure~\ref{fig:grid-coloring}. Reading the rows and columns gives exactly the six edges in \eqref{eq:grid-normal}.
\end{proof}

\FloatBarrier

We now use this representation to show that $H$ contains no $3\times3$ grid.

\begin{lemma}\label{lem:grid}
The hypergraph $H$ contains no $3\times3$ grid.
\end{lemma}

\begin{proof}
Suppose that $H$ contains a grid. By Lemma~\ref{lem:grid-normal}, its six edges can be written as $(X_i,Y_j,Z_{ij})$ for $i\ne j$, where $Z_{ij}=Z_{ji}$. Let $t$ be the least coordinate at which the three $X_i$ are not all equal or the three $Y_i$ are not all equal. Such a coordinate exists since the $X_i$ are distinct. Put $a_i=(X_i)_t$ and $b_i=(Y_i)_t$ for $1\le i\le3$.

By the choice of $t$, the cubic terms in coordinate $t$ of all six edges agree (and are absent if $t=1$). Comparing coordinate $t$ in $Z_{ij}=Z_{ji}$ gives $a_i+b_j=a_j+b_i$ for $i\ne j$. Hence there exists $D\in\F_p$ such that
\begin{equation}\label{eq:linear-at-t}
b_i=a_i+D\qquad(1\le i\le3).
\end{equation}
The $a_i$ are not all equal, since otherwise the $b_i$ would also be equal, contrary to the choice of $t$.

We first show that $t<d$. Otherwise, the condition $y_d=x_d$ gives $b_j=a_i$ for all $i\ne j$. Therefore, we have $a_1=a_2=a_3$ and $b_1=b_2=b_3$, contradicting the choice of $t$.

We next show that $a_1,a_2,a_3$ are pairwise distinct. Otherwise, after relabeling, we may assume that $a_1=a_2\ne a_3$. Since $t<d$, the three elements
\begin{equation*}
b_2-a_1=D,\qquad b_3-a_1=D+(a_3-a_1),\qquad b_1-a_3=D-(a_3-a_1)
\end{equation*}
all belong to $S$, contradicting its progression-free property.

Finally, comparing coordinate $t+1$ in $Z_{ij}=Z_{ji}$ for $(i,j)=(1,2),(2,3),(3,1)$ gives
\begin{equation*}
(X_i)_{t+1}+(Y_j)_{t+1}+a_i b_j^2
=(X_j)_{t+1}+(Y_i)_{t+1}+a_j b_i^2.
\end{equation*}
Summing these three equations, the linear terms cancel, and we obtain
\begin{equation*}
0=\begin{vmatrix}
1&1&1\\
a_1&a_2&a_3\\
b_1^2&b_2^2&b_3^2
\end{vmatrix}
=\begin{vmatrix}
1&1&1\\
a_1&a_2&a_3\\
a_1^2&a_2^2&a_3^2
\end{vmatrix}
=(a_2-a_1)(a_3-a_1)(a_3-a_2),
\end{equation*}
where the second equality follows from \eqref{eq:linear-at-t} by subtracting $2D$ times the second row and $D^2$ times the first row from the third row. This contradicts the pairwise distinctness of $a_1,a_2,a_3$.
\end{proof}

\noindent\emph{Remark.}
The cubic term in \eqref{eq:map} is essential to the preceding argument. If $x_{j-1}y_{j-1}^{\,2}$ were replaced by $x_{j-1}y_{j-1}$, the same calculation would yield
\begin{equation*}
0=\begin{vmatrix}
1&1&1\\
a_1&a_2&a_3\\
b_1&b_2&b_3
\end{vmatrix},
\end{equation*}
which holds identically by $b_i=a_i+D$. Thus this quadratic replacement does not yield a contradiction by the same argument, whereas the cubic term produces a nonzero Vandermonde determinant.

\begin{proof}[Proof of Theorem~\ref{thm:main}]
Let $n$ be sufficiently large, put $d=\lfloor(\log n)^{1/3}\rfloor$, and set $Q=(n/3)^{1/d}$. Since $Q\to\infty$ as $n\to\infty$, Bertrand's postulate gives a sufficiently large prime $p$ with $Q/2<p\le Q$. Choose a set $S\subseteq\F_p$ satisfying \eqref{eq:behrend}, and construct $H=H(p,d,S)$.

By Lemma~\ref{lem:linear-count}, $H$ is linear and has $3p^d\le n$ vertices. Lemmas~\ref{lem:triangle} and \ref{lem:grid} show that $H$ contains neither a loose triangle nor a grid. Hence, by Lemma~\ref{lem:classification}, $H$ is $(9,6)$-free. Adding isolated vertices gives an $n$-vertex linear tripartite $(9,6)$-free hypergraph with the same number of edges.

By Lemma~\ref{lem:linear-count} and \eqref{eq:behrend},
\begin{align*}
\log|E(H)|
&\ge (2d-1)\log p-C_0(d-1)\sqrt{\log p}\\
&\ge 2\log(n/3)-2d\log2-\frac{\log(n/3)}{d}-C_0\sqrt{d\log(n/3)}\\
&\ge 2\log n-O\bigl((\log n)^{2/3}\bigr),
\end{align*}
where the second inequality uses $Q/2<p\le Q$ and $(d-1)\sqrt{\log p}\le\sqrt{d\log(n/3)}$, and the last follows from $d=\Theta((\log n)^{1/3})$. The choice $d=\lfloor(\log n)^{1/3}\rfloor$ balances the two main losses $(\log n)/d$ and $\sqrt{d\log n}$, both of order $(\log n)^{2/3}$. Exponentiating gives the asserted bound, and in particular $f_3(n,9,6)\ge n^{2-o(1)}$.
\end{proof}

\section{Higher uniformities}\label{sec:higher}
The Vandermonde determinant argument used to exclude $3\times3$ grids can also be extended to higher uniformities. We combine it with a classification theorem due to Ge and Shangguan~\cite{GeShangguan2021} to prove Theorem~\ref{thm:higher}. Throughout this section we fix $r\ge4$. 

We first recall the relevant definition from \cite{GeShangguan2021}. A \emph{rainbow cycle} of length $\ell$ in an $r$-partite hypergraph consists of distinct edges $E_1,\ldots,E_\ell$ and vertices $w_i\in E_{i-1}\cap E_i$, with indices taken modulo $\ell$, such that the $w_i$ lie in pairwise distinct vertex parts.

\begin{lemma}[ {\cite[Theorem~5.6]{GeShangguan2021}}]\label{higher:lem:classification}
Let $F$ be a linear $r$-partite $r$-uniform hypergraph with no rainbow cycle of length three or four. If six edges of $F$ span at most $6r-9$ vertices, then they span exactly $6r-9$ vertices. Their projection onto three of the vertex parts is a $3\times3$ grid, and every vertex of their union outside these three parts belongs to exactly one of the six edges.
\end{lemma}

Once linearity is established, this classification reduces the proof of $(6r-9,6)$-freeness to excluding rainbow cycles of lengths three and four and $3\times3$ grids in projections onto three vertex parts.

\subsection{Solution-free sets and the construction}
As in \cite{GeShangguan2021}, we use solution-free sets to exclude rainbow cycles of lengths three and four. Let $R$ be a set of $r$ distinct elements of either $\mathbb Z$ or a field. We say that a set $S$ is \emph{$R_4$-solution-free} if, for $\ell\in\{3,4\}$ and pairwise distinct $c_1,\ldots,c_\ell\in R$, the equation
\begin{equation}\label{higher:eq:solution-free}
\sum_{i=1}^{\ell}(c_{i+1}-c_i)s_i=0,
\qquad c_{\ell+1}=c_1,\qquad s_1,\ldots,s_\ell\in S,
\end{equation}
implies $s_1=\cdots=s_\ell$. Write $\alpha_r=1-1/2^{r-1}$.

\begin{lemma}[Quantitative form of {\cite[Theorem~3.9]{GeShangguan2021}}]\label{higher:lem:integer-input}
There are constants $A_r,B_r>0$ such that, for every sufficiently large integer $N$, there exist an $r$-element set $R\subseteq\mathbb Z_{\ge0}$ and an $R_4$-solution-free set $S\subseteq\{1,\ldots,N\}$ satisfying
\begin{equation}\label{higher:eq:integer-input}
\max R\le\exp\bigl(A_r\sqrt{\log N}\bigr),
\qquad\textup{and}\qquad
|S|\ge N\exp\bigl(-B_r(\log N)^{\alpha_r}\bigr).
\end{equation}
\end{lemma}

We adapt this construction to the finite field $\F_p$. The idea is to choose the integer sets $R,S$ sufficiently small relative to $p$ so that reduction modulo $p$ preserves the solution-free property, without changing the asymptotic bound on $|S|$.

\begin{lemma}\label{higher:lem:field-input}
There is a constant $C_r>0$ such that, for every sufficiently large prime $p$, there exist an $r$-element set $R\subseteq\F_p$ and an $R_4$-solution-free set $S\subseteq\F_p$ with
\begin{equation}\label{higher:eq:field-input}
|S|\ge p\exp\bigl(-C_r(\log p)^{\alpha_r}\bigr).
\end{equation}
\end{lemma}

\begin{proof}
Choose a constant $K>A_r+1$ and put $N=\lfloor p\exp(-K\sqrt{\log p})\rfloor$. Since $N\to\infty$ as $p\to\infty$, Lemma~\ref{higher:lem:integer-input} gives sets $R\subseteq\mathbb Z_{\ge0}$ and $S\subseteq\{1,\ldots,N\}$ satisfying \eqref{higher:eq:integer-input}. For sufficiently large $p$, we have $4N\max R\le 4p\exp(-(K-A_r)\sqrt{\log p})<p$. In particular, reduction modulo $p$ is injective on $R$ and on $S$. Moreover, every sum in \eqref{higher:eq:solution-free}, regarded as an integer, has absolute value at most $4N\max R<p$. Hence, if such a sum vanishes modulo $p$, it also vanishes over the integers, forcing all the $s_i$ to be equal. Thus the reductions of $R$ and $S$ modulo $p$ satisfy the required $R_4$-solution-free property.

Since $N\ge (p/2)\exp(-K\sqrt{\log p})$ for sufficiently large $p$, Lemma~\ref{higher:lem:integer-input} gives
\begin{equation*}
|S|\ge \frac p2\exp\bigl(-K\sqrt{\log p}-B_r(\log p)^{\alpha_r}\bigr).
\end{equation*}
Since $\alpha_r>1/2$, this implies \eqref{higher:eq:field-input} for a suitable constant $C_r>0$.
\end{proof}

Below we present our construction for higher uniformities. Fix a prime $p>3$, an integer $d\ge2$, an $r$-element set $R\subseteq\F_p$, and an $R_4$-solution-free set $S\subseteq\F_p$. For each $b\in R$, let $V_b$ be a copy of $\F_p^d$. For $u,m\in\F_p^d$, define $v_b(u,m)\in V_b$ by
\begin{equation}\label{higher:eq:map}
\bigl(v_b(u,m)\bigr)_1=u_1+bm_1,
\qquad
\bigl(v_b(u,m)\bigr)_j=u_j+bm_j+b^2m_{j-1}^{\,3}
\quad(2\le j\le d).
\end{equation}
Let $H_r=H_r(p,d,R,S)$ be the $r$-partite $r$-uniform hypergraph with edges
\begin{equation}\label{higher:eq:edges}
E(u,m)=\{v_b(u,m):b\in R\},
\qquad
u\in\F_p^d,\quad
m\in S^{d-1}\times\{0\}.
\end{equation}
As in the three-uniform construction, the nonlinear term in coordinate $j$ depends on coordinate $j-1$. The condition $m_d=0$ ensures that, for any connected configuration with at least two edges, the first coordinate at which their parameter pairs $(u_j,m_j)$ are not all equal is less than $d$.

The next lemma establishes some basic properties of $H_r$.

\begin{lemma}\label{higher:lem:linear-count}
The hypergraph $H_r$ is linear, has $rp^d$ vertices, and has exactly $p^d|S|^{d-1}$ edges.
\end{lemma}

\begin{proof}
The number of vertices follows directly from the definition. There are $p^d|S|^{d-1}$ admissible pairs $(u,m)$. We show that any two vertices of an edge $E(u,m)$ uniquely determine its parameters $u,m$.

Fix distinct $a,b\in R$. The first coordinates of $v_a(u,m)$ and $v_b(u,m)$ are $u_1+am_1$ and $u_1+bm_1$, respectively, which uniquely determine $u_1,m_1$ since $a\ne b$. For $2\le j\le d$, suppose that $m_{j-1}$ has been determined. Subtracting $a^2m_{j-1}^3$ and $b^2m_{j-1}^3$ from the respective $j$th coordinates gives $u_j+am_j$ and $u_j+bm_j$, which uniquely determine $u_j,m_j$. Thus all parameters are recovered successively.

Consequently, distinct admissible pairs $(u,m)$ give distinct edges, and two distinct edges cannot share two vertices. This proves the asserted number of edges and linearity.
\end{proof}

\subsection{Excluding rainbow cycles and grids}

We now verify the two forbidden-configuration properties required by Lemma~\ref{higher:lem:classification}.

\begin{lemma}\label{higher:lem:cycles}
The hypergraph $H_r$ contains no rainbow cycle of length three or four.
\end{lemma}

\begin{proof}
Suppose that $E_i=E(u^{(i)},m^{(i)})$, $1\le i\le\ell$, form such a cycle, where $\ell\in\{3,4\}$ and $w_i\in V_{c_i}$ for distinct $c_1,\ldots,c_\ell\in R$. Let $t$ be the least coordinate for which the pairs $(u_t^{(i)},m_t^{(i)})$ are not all equal. Such a coordinate exists because the edges are distinct. Put $U_i=u_t^{(i)}$ and $M_i=m_t^{(i)}$.

At each intersection $w_i$, the nonlinear terms in coordinate $t$ agree by the minimality of $t$ (and are absent if $t=1$). Hence $U_{i-1}+c_iM_{i-1}=U_i+c_iM_i$ for every $i$. Summing these equalities gives
\begin{equation*}
\sum_{i=1}^{\ell}(c_{i+1}-c_i)M_i=0.
\end{equation*}
If $t<d$, all $M_i$ belong to $S$, so \eqref{higher:eq:solution-free} implies that they are equal. If $t=d$, they are all zero by \eqref{higher:eq:edges}. In either case, the intersection equalities imply that all $U_i$ are equal, contradicting the choice of $t$.
\end{proof}

We next exclude any grid in the projection onto any three vertex parts. For distinct $a,b,c\in R$, this projection has edges $E\cap(V_a\cup V_b\cup V_c)$, where $E$ ranges over the edges of $H_r$. By linearity, these projected triples are distinct and form a linear tripartite hypergraph.

\begin{lemma}\label{higher:lem:grid}
For any three distinct $a,b,c\in R$, the projection of $H_r$ onto $V_a\cup V_b\cup V_c$ contains no $3\times3$ grid.
\end{lemma}

\begin{proof}
Suppose that a grid occurs in this projection. By Lemma~\ref{lem:grid-normal}, its six triples can be written as $(X_i,Y_j,Z_{ij})$ for $1\le i,j\le3$ with $i\ne j$, where $Z_{ji}=Z_{ij}$, $X_i\in V_a$, $Y_i\in V_b$, and $Z_{ij}\in V_c$. Denote the corresponding edges of $H_r$ by $E_{ij}=E(u^{(ij)},m^{(ij)})$.

Let $t$ be the least coordinate for which the six pairs $(u_t^{(ij)},m_t^{(ij)})$ are not all equal. Such a coordinate exists because the six edges are distinct. We first show that $t<d$. Otherwise, $m_d^{(ij)}=0$ for every edge, and by the minimality of $t$, all six edges have the same value of $m_{d-1}$. Thus two edges sharing a vertex must have the same value of $u_d^{(ij)}$. Since the grid is connected, these six values are all equal, contradicting the choice of $t$.

By the minimality of $t$, the nonlinear terms in coordinate $t$ agree across all six edges within each vertex part (and are absent if $t=1$). Subtract these common terms from the corresponding coordinates, and denote the resulting values for $X_i,Y_i,Z_{ij}$ by $A_i,B_i,C_{ij}$, respectively. Then, for every $i\ne j$,
\begin{equation*}
A_i=u_t^{(ij)}+am_t^{(ij)},\qquad
B_j=u_t^{(ij)}+bm_t^{(ij)},\qquad
C_{ij}=u_t^{(ij)}+cm_t^{(ij)}.
\end{equation*}
By \eqref{higher:eq:map}, we have $(b-a)C_{ij}=(b-c)A_i+(c-a)B_j$. Since $C_{ij}=C_{ji}$, there is a $D\in\F_p$ such that
\begin{equation}\label{higher:eq:grid-linear}
B_i=qA_i+D\quad(1\le i\le3),
\qquad q=\frac{b-c}{c-a}.
\end{equation}
We therefore have
\begin{equation}\label{higher:eq:grid-slopes}
m_t^{(ij)}=\frac{B_j-A_i}{b-a}
=\frac{qA_j-A_i+D}{b-a}\in S
\qquad(i\ne j).
\end{equation}
The $A_i$ are not all equal. Otherwise, the $B_i$ would also all be equal. Since $a\ne b$, the relations $A_i=u_t^{(ij)}+am_t^{(ij)}$ and $B_j=u_t^{(ij)}+bm_t^{(ij)}$ would then force all six pairs $(u_t^{(ij)},m_t^{(ij)})$ to be equal, contradicting the choice of $t$.

In fact, $A_1,A_2,A_3$ are pairwise distinct. Otherwise, after relabeling, we may assume $A_1=A_2\ne A_3$. By \eqref{higher:eq:grid-linear}, we have $B_1=B_2$. Equation~\eqref{higher:eq:grid-slopes} then gives
\begin{equation*}
m_t^{(13)}-m_t^{(12)}=\frac{q(A_3-A_1)}{b-a},
\qquad
m_t^{(12)}-m_t^{(31)}=\frac{A_3-A_1}{b-a}.
\end{equation*}
Since $q=(b-c)/(c-a)$, these equalities yield $(c-a)m_t^{(13)}+(b-c)m_t^{(31)}+(a-b)m_t^{(12)}=0$. Applying \eqref{higher:eq:solution-free} with $\ell=3$ and $(c_1,c_2,c_3)=(a,c,b)$, we obtain $m_t^{(13)}=m_t^{(31)}=m_t^{(12)}$, contradicting $A_1\ne A_3$.

We now compare coordinate $t+1$. For every edge $E(u,m)$, the definition gives
\begin{equation}\label{higher:eq:interpolation}
(b-a)\bigl(v_c(u,m)\bigr)_{t+1}
-(b-c)\bigl(v_a(u,m)\bigr)_{t+1}
-(c-a)\bigl(v_b(u,m)\bigr)_{t+1}
=(b-a)(c-a)(c-b)m_t^3.
\end{equation}
Apply \eqref{higher:eq:interpolation} to each of the six edges. Subtract the equation for $E_{21}$ from that for $E_{12}$, and similarly for the pairs $(E_{23},E_{32})$ and $(E_{31},E_{13})$. Adding the three resulting equations, all terms involving vertex coordinates cancel, since $Z_{ij}=Z_{ji}$ and each $X_i$ and $Y_i$ occurs once with each sign. Since $(b-a)(c-a)(c-b)\ne0$, we obtain
\begin{equation}\label{higher:eq:cubic-sum}
0=\bigl((m_t^{(12)})^3-(m_t^{(21)})^3\bigr)
+\bigl((m_t^{(23)})^3-(m_t^{(32)})^3\bigr)
+\bigl((m_t^{(31)})^3-(m_t^{(13)})^3\bigr).
\end{equation}
Substituting \eqref{higher:eq:grid-slopes} into \eqref{higher:eq:cubic-sum}, the terms involving $D$ and the pure cubic terms cancel in the cyclic sum. The remaining terms give
\begin{align*}
0&=-\frac{3q(q+1)}{(b-a)^3}
\bigl(A_1A_2(A_2-A_1)+A_2A_3(A_3-A_2)+A_3A_1(A_1-A_3)\bigr)\\
&=-\frac{3q(q+1)}{(b-a)^3}
\begin{vmatrix}
1&1&1\\
A_1&A_2&A_3\\
A_1^2&A_2^2&A_3^2
\end{vmatrix}.
\end{align*}
Since $p>3$, $q\ne0$, $q+1=(b-a)/(c-a)\ne0$, and the $A_i$ are pairwise distinct, the right-hand side is nonzero, a contradiction.
\end{proof}

We now combine the preceding lemmas to establish the required $(6r-9,6)$-free property.

\begin{proposition}\label{higher:prop:avoidance}
The hypergraph $H_r$ is $(6r-9,6)$-free.
\end{proposition}

\begin{proof}
By Lemmas~\ref{higher:lem:linear-count} and~\ref{higher:lem:cycles}, $H_r$ satisfies the hypotheses of Lemma~\ref{higher:lem:classification}. Any six edges spanning at most $6r-9$ vertices would therefore give a grid in the projection onto three parts, contrary to Lemma~\ref{higher:lem:grid}.
\end{proof}

\begin{proof}[Proof of Theorem~\ref{thm:higher}]
Let $n$ be sufficiently large, put $L=\log(n/r)$ and $\beta_r=(1-\alpha_r)/(2-\alpha_r)=1/(2^{r-1}+1)$, and choose $d=\lfloor L^{\beta_r}\rfloor$ and $Q=(n/r)^{1/d}$. Since $r$ is fixed, we have $d\ge2$ and $Q\to\infty$ as $n\to\infty$. By Bertrand's postulate, there is a sufficiently large prime $p$ with $Q/2<p\le Q$. Choose $R,S\subseteq\F_p$ as in Lemma~\ref{higher:lem:field-input}, and construct $H_r$.

By Lemma~\ref{higher:lem:linear-count} and Proposition~\ref{higher:prop:avoidance}, the hypergraph $H_r$ is linear, $r$-partite, and $(6r-9,6)$-free, with $rp^d\le n$ vertices. Moreover,
\begin{align*}
\log|E(H_r)|
&\ge(2d-1)\log p-C_r(d-1)(\log p)^{\alpha_r}\\
&\ge2L-\frac Ld-2d\log2-C_rL^{\alpha_r}d^{1-\alpha_r}\\
&\ge2\log n-O_r\bigl((\log n)^{\theta_r}\bigr).
\end{align*}
The second inequality follows from $L/d-\log2<\log p\le L/d$. For the last inequality, note that $d=\Theta_r(L^{\beta_r})$ and $1-\beta_r=\alpha_r+(1-\alpha_r)\beta_r=\theta_r$, while $\beta_r<\theta_r$. Exponentiating and adding isolated vertices gives the required hypergraph.
\end{proof}

\section*{Acknowledgements}
Chong Shangguan is supported by the National Natural Science Foundation of China under Grant Nos.~12571352 and 12231014, and the Fundamental Research Funds for the Central Universities. He thanks Baoyan Feng, Jiaxi Nie, and Yulin Yang for helpful discussions of the papers~\cite{GishbolinerShapira2022,GowersJanzer2021,Solymosi2026}.

During the preparation of this manuscript, the authors used OpenAI's ChatGPT (GPT-6) as an auxiliary tool for exploratory mathematical discussions, checking calculations, and improving the exposition. The authors take full responsibility for the manuscript.

\appendix
\section{Further forbidden configurations}\label{app:further}
We continue with $H=H(p,d,S)$ from Section~\ref{sec:construction}. The same construction is $(e+3,e)$-free for every $3\le e\le11$. The cases $e\le8$ follow from the absence of loose triangles and grids together with elementary degree arguments (see Proposition~\ref{further:through-eight}). For $e=9$, the additive relations at the first differing coordinate give a contradiction (see Proposition~\ref{further:nine}). The case $e=10$ requires a further identity obtained by taking a weighted sum of the cubic terms (see Lemma~\ref{further:first-coordinate} and Proposition~\ref{further:ten}), while the case $e=11$ follows from a counting argument (see Proposition~\ref{further:eleven}).

All degrees below are taken within the configuration under consideration.

\begin{lemma}\label{further:degree}
Let $F$ be a linear $3$-uniform hypergraph with $m$ edges, no loose triangle, and minimum degree at least two. Then every vertex of $F$ has degree at most $\lfloor m/3\rfloor$.
\end{lemma}

\begin{proof}
Let $v$ have degree $q$, and let $E_1,\ldots,E_q$ be its incident edges. By linearity, these edges contain $2q$ distinct vertices other than $v$. Moreover, every edge outside $\{E_1,\ldots,E_q\}$ meets their union in at most one vertex. Indeed, two intersections with the same $E_i$ violate linearity, while intersections with two different $E_i$ form a loose triangle. Each of the $2q$ vertices other than $v$ has degree at least two, so at least $2q$ edges outside the star are required. Hence $m-q\ge2q$.
\end{proof}

\begin{lemma}\label{further:reduction}
Let $m\ge4$, and let $G$ be a linear $3$-uniform hypergraph with no loose triangle. Suppose that $G$ is $(m+2,m-1)$-free. If $m$ edges of $G$ span at most $m+3$ vertices, then they span exactly $m+3$ vertices, and their configuration has minimum degree at least two and maximum degree at most $\lfloor m/3\rfloor$.
\end{lemma}

\begin{proof}
If their union had at most $m+2$ vertices, deleting any edge would contradict the $(m+2,m-1)$-freeness of $G$. If a vertex had degree one, deleting its unique incident edge would again leave $m-1$ edges on at most $m+2$ vertices. The maximum-degree bound now follows from Lemma~\ref{further:degree}.
\end{proof}

\begin{proposition}\label{further:through-eight}
The hypergraph $H$ is $(e+3,e)$-free for every $3\le e\le8$.
\end{proposition}

\begin{proof}
Three edges of a linear $3$-uniform hypergraph span at most six vertices only if they form a loose triangle. Thus Lemma~\ref{lem:triangle} gives the case $e=3$. Applying Lemma~\ref{further:reduction} successively with $m=4,5$ would give a forbidden configuration of minimum degree at least two and maximum degree at most one. This proves the cases $e=4,5$.

The case $e=6$ is the $(9,6)$-free property established in Section~\ref{sec:construction}. For $m=7$, and then for $m=8$, Lemma~\ref{further:reduction} would force every degree of a forbidden configuration to equal two. Writing $v$ for its number of vertices, we would have $3m=2v\le2(m+3)$, which is impossible for $m>6$.
\end{proof}

For the remaining cases, we examine the first coordinate at which the vertex labels differ. We say a configuration in $H$ is connected if its incidence graph is connected. The following lemma gives the resulting additive relations and extends the cancellation used in Lemma~\ref{lem:grid} to weighted sums over the edges of a configuration.


\begin{lemma}\label{further:first-coordinate}
Let $F$ be a connected configuration of at least two edges in $H$. Let $k$ be the least coordinate at which the vertices of $F$ in $V_1$, or those in $V_2$, do not all have the same coordinate value. Then $k<d$. At coordinate $k$, after subtracting a common constant from the labels in $V_3$, every edge satisfies
\begin{equation}\label{further:additive-layer}
z=x+y,\qquad y-x\in S.
\end{equation}
Moreover, if weights $\lambda_E\in\F_p$, $E\in E(F)$, satisfy
\begin{equation}\label{further:vertex-balance}
\sum_{E\ni v}\lambda_E=0\qquad\text{for every }v\in V(F),
\end{equation}
then the coordinate-$k$ labels in the original first two parts satisfy
\begin{equation}\label{further:weighted-cubic}
\sum_{E=(x,y,z)\in E(F)}\lambda_E x y^2=0.
\end{equation}
Equivalently,
\begin{equation}\label{further:slope-cubic}
\sum_{E=(x,y,z)\in E(F)}\lambda_E (y-x)^3=0.
\end{equation}
\end{lemma}

\begin{proof}
Such a coordinate exists: otherwise all edges would have the same vertices in $V_1$ and $V_2$, and hence would be the same edge. At every coordinate before $k$, all $V_1$ labels agree and all $V_2$ labels agree. Consequently, the cubic terms in coordinate $k$ of \eqref{eq:map} have a common value $C$; we set $C=0$ if $k=1$. To distinguish vertices from their scalar labels, write an edge as $E=(X_E,Y_E,Z_E)$, where $X_E,Y_E,Z_E$ are vectors in the respective vertex parts.

Suppose that $k=d$. By \eqref{eq:admissible}, each edge satisfies $(X_E)_d=(Y_E)_d$. Then, by \eqref{eq:map}, we have $(Z_E)_d=2(X_E)_d+C$. If two edges share a vertex in $V_1$ or $V_2$, their values of $(X_E)_d=(Y_E)_d$ are equal. If they share a vertex in $V_3$, the shared vertex in $V_3$ implies that $2(X_E)_d+C=2(X_{E'})_d+C$. Since $p>2$, this again implies $(X_E)_d=(X_{E'})_d$. Connectedness therefore forces this value to be the same on every edge, contrary to the definition of $k$. Thus $k<d$. Set
\[
x_E=(X_E)_k,\qquad y_E=(Y_E)_k,\qquad z_E=(Z_E)_k-C.
\]
Since $k<d$, equation~\eqref{eq:map} and the admissibility
conditions~\eqref{eq:admissible} give \eqref{further:additive-layer}.

Note that by \eqref{further:vertex-balance}, for each $i\in\{1,2,3\}$ and any function $g:V(F)\cap V_i\to\F_p$, we have
\begin{equation*}
\sum_{E\in E(F)}\lambda_E\sum_{v\in E\cap V_i}g(v)
=\sum_{v\in V(F)\cap V_i}g(v)\left(\sum_{E\ni v}\lambda_E\right)=0.
\end{equation*}
Hence, multiplying the equation in \eqref{eq:map} for coordinate-$(k+1)$
\[
(Z_E)_{k+1}=(X_E)_{k+1}+(Y_E)_{k+1}+x_Ey_E^2.
\]
by $\lambda_E$ and summing over $E\in E(F)$, the linear terms vanish separately. This results in $\sum_{E\in E(F)}\lambda_Ex_Ey_E^2=0$, which is \eqref{further:weighted-cubic}.

Finally, since $z_E=x_E+y_E$, we have $6x_Ey_E^2=z_E^3-(y_E-x_E)^3-2x_E^3$. By the same argument, the weighted sums of $x_E^3$ and $z_E^3$ also vanish. Since $p>3$, this proves that \eqref{further:weighted-cubic} and \eqref{further:slope-cubic} are equivalent.
\end{proof}


For nine edges, the additive relations alone give a contradiction. The required configuration is derived directly in the proof.

\begin{proposition}\label{further:nine}
The hypergraph $H$ is $(12,9)$-free.
\end{proposition}

\begin{proof}
Suppose that nine edges of $H$ span at most twelve vertices, and let $F$ be their configuration. By Proposition~\ref{further:through-eight} and Lemma~\ref{further:reduction}, $F$ has twelve vertices and minimum degree at least two. Write $V_i(F)=V(F)\cap V_i$. Since every edge meets each vertex part in exactly one vertex,
\[
\sum_{v\in V_i(F)}\deg_F(v)=9.
\]
Thus, $2|V_i(F)|\leq 9$, which implies that each part contains at most four vertices. As $F$ has twelve vertices in total, each part contains exactly four vertices with degree sequence $(3,2,2,2)$.

Let $c$ be the vertex of degree three in $V_3(F)$. By linearity, the three edges through $c$ have distinct vertices in each of $V_1,V_2$. Write $V_1(F)=\{X_1,X_2,X_3,X_4\}$, $V_2(F)=\{Y_1,Y_2,Y_3,Y_4\}$, and label these edges as $(X_i,Y_i,c)$ for $1\le i\le3$.

Consider the $4\times4$ array with rows indexed by $V_1(F)$ and columns indexed by $V_2(F)$. We call a cell $(i,j)$ occupied if $(X_i,Y_j,Z)\in E(F)$ for some $Z\in V_3(F)$. By linearity, each occupied cell represents a unique edge. Then, no cell $(i,j)$ with $1\le i,j\le3$ and $i\ne j$ is occupied: an edge in this cell would meet the two diagonal edges $(X_i,Y_i,c)$ and $(X_j,Y_j,c)$ at $X_i$ and $Y_j$, respectively, giving a loose triangle.

Now, we consider the other six edges. They must therefore occupy six of the seven cells in the fourth row or fourth column. If the corner $(4,4)$ were occupied, it would be counted twice in the sum of the degrees of $X_4$ and $Y_4$, while each of the other five edges would be counted once. This sum would be $2+5=7$, forcing one of $X_4$ and $Y_4$ to have degree at least four, contrary to the degree sequences above. Hence the corner is empty and all six non-corner cells are occupied.

Each of the three remaining vertices of $V_3(F)$ has degree two. We call the two occupied cells representing the edges incident with such a vertex a pair. These three pairs partition the six non-corner cells. By linearity, the two cells of a pair cannot lie in the same row or column. Consequently, each pair consists of one cell of the fourth row and one cell of the fourth column. If $(4,i)$ were paired with $(i,4)$, with common third-part vertex $e$, the edges $(X_4,Y_i,e), ~(X_i,Y_4,e),~(X_i,Y_i,c)$ would meet pairwise at the three distinct vertices $e,X_i,Y_i$, forming a loose triangle.

The pairing therefore defines a permutation $\sigma$ of $\{1,2,3\}$, where $(4,i)$ is paired with $(\sigma(i),4)$, and $\sigma(i)\ne i$ for every $i$. Such a permutation is a $3$-cycle. After simultaneously relabeling the first three rows and columns, we may assume that the pairs
are
\[
\{(4,1),(3,4)\},\qquad
\{(4,2),(1,4)\},\qquad
\{(4,3),(2,4)\}.
\]
Put $a=X_4$, $b=Y_4$, $d_i=Y_i$, and
$(f_1,f_2,f_3)=(X_3,X_1,X_2)$. Read subscripts on $d_i,e_i,f_i$ modulo three. Let $e_i$ be the common third-part vertex of the pair containing $(4,i)$. The fourth row and fourth column then give the first two families below, while the original diagonal cells give the third, since $X_i=f_{i+1}$:
\begin{equation}\label{further:nine-configuration}
(a,d_i,e_i),\qquad(f_i,b,e_i),\qquad(f_{i+1},d_i,c)
\qquad(1\le i\le3).
\end{equation}
This configuration is connected: the three edges containing $c$ lie in one component, and each of the other six edges meets one of them at some $d_i$ or $f_i$.

Now, we apply Lemma~\ref{further:first-coordinate} to $F$, and use the same letters for the scalar labels in \eqref{further:additive-layer}. The first two families of edges give $e_i=a+d_i=f_i+b$, and hence $f_i=a+d_i-b$. The third family then gives
\[
d_i+d_{i+1}=c-a+b\qquad(1\le i\le3).
\]
Subtracting the equalities for $i=1,2$ gives $d_1=d_3$, and subtracting those for $i=2,3$ gives $d_2=d_1$. Thus $d_1=d_2=d_3=D$, and hence $f_1=f_2=f_3=f=a+D-b$. Put $s=D-a$ and $\eta=D-b=f-a$. The three types of edges in \eqref{further:nine-configuration} have differences
\[
D-a=s,\qquad b-f=s-2\eta,\qquad D-f=s-\eta.
\]
All three belong to $S$. If $\eta\ne0$, they form a nontrivial three-term arithmetic progression, contradicting the progression-freeness of $S$. Therefore $\eta=0$, and hence $f=a$ and $D=b$, making all scalar labels constant within each of $V_1(F),V_2(F)$, contrary to the choice of the first differing coordinate.
\end{proof}

For ten edges, we first describe the possible configuration. 

\begin{lemma}\label{further:ten-configuration}
Let $F$ be a linear tripartite $3$-uniform hypergraph with ten edges and thirteen vertices, whose degrees are all two or three. If $F$ contains neither a loose triangle nor a $3\times3$ grid, then its vertex parts have sizes $4,4,5$. After naming the two parts containing four vertices $U=\{U_0,U_1,U_2,U_3\}$ and $V=\{V_0,V_1,V_2,V_3\}$, its ten edges are represented by the following $U$-by-$V$ array, whose nonempty entry is the vertex in the third part $W$:
\begin{equation}\label{further:ten-array}
\begin{pmatrix}
W_0&W_1&W_2&\cdot\\
W_3&W_4&\cdot&W_2\\
\cdot&\cdot&W_4&W_0\\
\cdot&\cdot&W_3&W_1
\end{pmatrix}.
\end{equation}
In particular, $F$ is connected.
\end{lemma}

\begin{proof}
The degree sum in each part is ten, and every degree is two or three. Thus each part has four or five vertices, and the total of thirteen forces part sizes $4,4,5$. The degree sequences in the two four-vertex parts are $(3,3,2,2)$, while every vertex of $W$ has degree two.

As in the proof of Proposition~\ref{further:nine}, we represent edges by occupied cells of a $4\times4$ array. Each vertex $w\in W$ corresponds to a pair of cells $(i,j),(i',j')$ lying in distinct rows and distinct columns. The other two corners of this rectangle, $(i,j')$ and $(i',j)$, referred to as the cross cells, must both be empty. Otherwise, an edge in $(i,j')$ would have third-part vertex different from $w$ by linearity, and would form a loose triangle with the paired edges. The other cross cell is excluded in the same way.

A row of degree three has a unique empty cell, say in column $q$. If an occupied cell in this row is paired with a cell in column $q'$, the cross cell in the original row and column $q'$ must be empty. Thus $q'=q$, so all three partners of the occupied cells in this row lie in column $q$. The partners are distinct, and column $q$ already has an empty cell in the original row; hence its degree is exactly three. The two rows of degree three cannot miss the same column, since that column would then have two empty cells and degree at most two, contradicting the preceding conclusion that it has degree three.

Relabel the degree-three rows as $0,1$, and the columns they miss as $3,2$, respectively. The first two rows then have the form $(*,*,*,\cdot)$ and $(*,*,\cdot,*)$, where we use $*$ to denote an occupied cell and $\cdot$ to denote an empty cell. Columns $2,3$ both have degree three and each currently has one occupied cell, so all four cells in rows $2,3$ and columns $2,3$ must be occupied. Since rows $2,3$ have degree two, their remaining cells are empty. Thus the occupied cells have the form
\[
\begin{pmatrix}
*&*&*&\cdot\\
*&*&\cdot&*\\
\cdot&\cdot&*&*\\
\cdot&\cdot&*&*
\end{pmatrix}.
\]

The partner of $(0,2)$ lies in column $3$. It cannot be $(2,3)$ or $(3,3)$, because the corresponding cross cell $(2,2)$ or $(3,2)$ is occupied. Hence $(0,2)$ pairs with $(1,3)$. The cells $(0,0),(0,1)$ therefore pair bijectively with $(2,3),(3,3)$, and $(1,0),(1,1)$ pair bijectively with $(2,2),(3,2)$.

Suppose that these two bijections assign a common column $j\in\{0,1\}$ to the same row $i\in\{2,3\}$. Then, the three pairs are $\{(0,2),(1,3)\}$, $\{(0,j),(i,3)\}$, $\{(1,j),(i,2)\}$. Let their third-part vertices be $w_0,w_1,w_2$, respectively. These vertices are distinct, since the three pairs are disjoint and every vertex of $W$ has degree two. Partition the six
corresponding edges into
\begin{align*}
R_1=(U_0,V_2,w_0),&~C_1=(U_1,V_3,w_0),\\
R_2=(U_i,V_3,w_1),&~C_2=(U_0,V_j,w_1),\\
R_3=(U_1,V_j,w_2),&~C_3=(U_i,V_2,w_2).
\end{align*}
The edges $R_1,R_2,R_3$ are pairwise disjoint, as are $C_1,C_2,C_3$.
Their intersections are given by the following array, with each entry
denoting the unique intersection vertex:
\[
\begin{array}{c|ccc}
 & C_1&C_2&C_3\\ \hline
R_1&w_0&U_0&V_2\\
R_2&V_3&w_1&U_i\\
R_3&U_1&V_j&w_2
\end{array}.
\]
All nine entries are distinct: $U_0,U_1,U_i$ are distinct, as are
$V_j,V_2,V_3$ and $w_0,w_1,w_2$. These six edges therefore form a
$3\times3$ grid, a contradiction.

Thus, the two bijections disagree at both columns. Interchanging rows $2$ and $3$ if necessary, we may assume that $(0,0)$ pairs with $(2,3)$ and $(0,1)$ with $(3,3)$. Then, the other bijection pairs $(1,0)$ with $(3,2)$ and $(1,1)$ with $(2,2)$. Label the vertices of $W$ corresponding to the pairs containing $(0,0),(0,1),(0,2),(1,0),(1,1)$ by $W_0,W_1,W_2,W_3,W_4$, respectively. This gives \eqref{further:ten-array}. Finally, the connectedness of $F$ follows directly from the form of the ten edges shown in array \eqref{further:ten-array}.
\end{proof}

\begin{proposition}\label{further:ten}
The hypergraph $H$ is $(13,10)$-free.
\end{proposition}

\begin{proof}
Suppose that ten edges of $H$ span at most thirteen vertices. Proposition~\ref{further:nine} and Lemma~\ref{further:reduction} imply that their configuration $F$ has thirteen vertices and all degrees equal to two or three. By Lemmas~\ref{lem:triangle}, \ref{lem:grid}, and~\ref{further:ten-configuration}, its edges have the form \eqref{further:ten-array}, and $F$ is connected.

Apply Lemma~\ref{further:first-coordinate}, and let $k$ be the first differing coordinate. For each edge $E=(X,Y,Z)$ with $X\in V_1$, $Y\in V_2$ and $Z\in V_3$, write $x=X_k,~y=Y_k,~z=Z_k-C$, where $C$ is the common constant from that lemma. These scalar labels lie in $\F_p$ and satisfy $z=x+y$ and $y-x\in S$. Write $V_i(F)=V(F)\cap V_i$. 

Denote the two four-vertex parts by $U,V$ and the five-vertex part by $W$ as in Lemma \ref{further:ten-configuration}. We order these parts and assign auxiliary scalar labels $u,v,w$
to their vertices according to the following table:
\[
\begin{array}{c|c|c}
\text{five-vertex part}
&(U, V,W)&(u,v,w)\\ \hline
V_3(F)&(V_1(F),V_2(F),V_3(F))&(x,y,z)\\
V_1(F)&(V_2(F),V_3(F),V_1(F))&(-y,z,x)\\
V_2(F)&(V_1(F),V_3(F),V_2(F))&(-x,z,y)
\end{array}.
\]
Thus, by $z=x+y$, every edge satisfies $w=u+v$.

Now, we label the vertices in these ordered parts as $U_i,V_j,W_\ell$ according to \eqref{further:ten-array}, and denote their auxiliary scalar labels by $u_i,v_j,w_\ell$, respectively. Then, an occupied cell $(i,j)$ carrying $W_\ell$ satisfies $u_i+v_j=w_\ell$. Hence, the five pairs of cells carrying the same vertex of $W$ give
\begin{align*}
u_0+v_0=u_2+v_3,&~u_0+v_1=u_3+v_3,\\
u_0+v_2=u_1+v_3,&~u_1+v_0=u_3+v_2,\\
u_1+v_1=u_2+v_2.
\end{align*}
Set $A=u_0$, $B=v_3$, and $t=v_2-v_3$. The first three relations yield $u_2=A+v_0-B$, $u_3=A+v_1-B$, $u_1=A+t$. Substituting these expressions into the fourth relation yields $v_0=v_1$; the fifth relation then holds as well. With $\delta=v_0-B$, we obtain
\begin{equation}\label{further:ten-values}
\begin{split}
(u_0,u_1,u_2,u_3)&=(A,A+t,A+\delta,A+\delta),\\
(v_0,v_1,v_2,v_3)&=(B+\delta,B+\delta,B+t,B).
\end{split}
\end{equation}
Here $(\delta,t)\ne(0,0)$. Otherwise, all $u_i$ would equal $A$, all $v_j$ would equal $B$, and all $w_\ell$ would equal $A+B$. Then the original labels at coordinate $k$ would be constant within each vertex part, contradicting the choice of $k$.

Next, give the nonempty cells of \eqref{further:ten-array} the weights
\[
\Lambda=
\begin{pmatrix}
1&1&-2&0\\
-1&-1&0&2\\
0&0&1&-1\\
0&0&1&-1
\end{pmatrix},
\]
where empty cells have weight zero. For each edge $E=(U_i,V_j,W_\ell)$, set $\lambda_E=\Lambda_{ij}$, with $0\le i,j\le3$. Clearly, each row and each column has weight sum zero. The two cells carrying $W_2$ have weights $-2,2$, and the two cells carrying any other fixed $W_\ell$ have weights $1,-1$. Thus the weights also sum to zero at every vertex of $W$.

For a polynomial $Q\in\F_p[u,v,w]$, define
\[
L(Q)=\sum_{E=(U_i,V_j,W_\ell)\in E(F)}
\lambda_E Q(u_i,v_j,w_\ell).
\]
If $Q$ depends only on $u$, say $Q=g(u)$, grouping by rows gives
\[
L(g(u))=\sum_{i=0}^3g(u_i)\left(\sum_{j=0}^3\Lambda_{ij}\right)=0.
\]
Grouping by columns or by the pairs carrying a fixed $W_\ell$ gives the same conclusion when $Q$ depends only on $v$ or only on $w$. Grouping the terms of $L(uv^2)$ by rows gives
\begin{align*}
L(uv^2) & = \sum_{0\leq i,j\leq 3}\Lambda_{ij}u_iv_j^2 = u_0(v_0^2+v_1^2-2v_2^2)+u_1(-v_0^2-v_1^2+2v_3^2)
 +(u_2+u_3)(v_2^2-v_3^2) \\
 & = 2A\bigl((B+\delta)^2-(B+t)^2\bigr)+2(A+t)\bigl(B^2-(B+\delta)^2\bigr)+2(A+\delta)\bigl((B+t)^2-B^2\bigr) \\
 & = 2\delta t(t-\delta). 
\end{align*}
Here we first group the terms by rows and then substitute \eqref{further:ten-values}; the final equality follows by expanding and canceling the terms involving $A$ and $B$. Then, by $L(u^3)=L(v^3)=L(w^3)=0$ and $w=u+v$ on every edge, we have
\[
0=L((u+v)^3)=3L(u^2v)+3L(uv^2).
\]
Since $p>3$, it follows that $L(u^2v)=-L(uv^2)$.

In the three choices of auxiliary labels above, the original difference
$y-x$ equals, respectively,
\[
v-u,\qquad -u-w=-2u-v,\qquad w+u=2u+v.
\]
Thus $y-x=\alpha u+\beta v$, where $(\alpha,\beta)\in\{(-1,1),(-2,-1),(2,1)\}$. The weights satisfy \eqref{further:vertex-balance} at the original vertices as well, so Lemma~\ref{further:first-coordinate} gives $L((y-x)^3)=0$. Then, using the vanishing pure-cube sums and $L(u^2v)=-L(uv^2)$, we obtain
\begin{equation*}
0=L((\alpha u+\beta v)^3)=3\alpha^2\beta L(u^2v)+3\alpha\beta^2L(uv^2)=6\alpha\beta(\beta-\alpha)\delta t(t-\delta).
\end{equation*}
For the three pairs $(\alpha,\beta)$, the coefficient $\alpha\beta(\beta-\alpha)$ is $-2,2,-2$, respectively. Since $p>3$, we conclude that $\delta t(t-\delta)=0$.

Recall that the original differences $y-x$ belong to $S$. Consequently, after subtracting $\alpha A+\beta B$, their values belong to $S-(\alpha A+\beta B)$, which is also progression-free. On the ten edges, these translated differences are given by the following six expressions:
\[
\beta\delta,\quad\alpha\delta,\quad\beta t,\quad\alpha t,
\quad\alpha t+\beta\delta,\quad\alpha\delta+\beta t.
\]
If $\delta=0$ or $t=0$, the other parameter is some $s\ne0$, and these values include $0,\alpha s,\beta s$. If neither parameter is zero, then $\delta t(t-\delta)=0$ gives $\delta=t=s\ne0$, and the values include $\alpha s,\beta s,(\alpha+\beta)s$. Since $s\ne0$ and $p>3$, for each choice of $(\alpha,\beta)\in\{(-1,1),(-2,-1),(2,1)\}$, each triple consists of three distinct terms of an arithmetic progression, a contradiction. This proves the proposition.
\end{proof}

The next case follows from a short counting argument.

\begin{proposition}\label{further:eleven}
The hypergraph $H$ is $(14,11)$-free.
\end{proposition}

\begin{proof}
Suppose that eleven edges of $H$ span at most fourteen vertices. By Proposition~\ref{further:ten} and Lemma~\ref{further:reduction}, their configuration has fourteen vertices and all degrees are two or three. The degree sum in each part is eleven, so the part sizes are $4,5,5$. Call the part containing four vertices $U$ and the other parts $V,W$. Their degree sequences are $(3,3,3,2)$, $(3,2,2,2,2)$, and $(3,2,2,2,2)$, respectively.

In the $4\times5$ $U$-by-$V$ array, label the three edges through the vertex of degree three in $W$ as $(0,0),(1,1),(2,2)$. Since there is no loose triangle, the six off-diagonal cells of this $3\times3$ block are empty. Each of columns $0,1,2$ has degree at least two, so row $3$ must contain all three cells $(3,0),(3,1),(3,2)$. Its degree is therefore three, and its other two cells are empty. The remaining five edges occupy five of the six cells in the $3\times2$ block in the upper right corner.

Each of the other four vertices of $W$ pairs two cells outside the initial diagonal. As in the proof of Lemma~\ref{further:ten-configuration}, the cells of a pair lie in distinct rows and columns and have both cross cells empty. Thus a pair cannot contain two of $(3,0),(3,1),(3,2)$, which lie in one row. Nor can it contain two cells of this $3\times2$ block, since that block has only one empty cell. Every pair must therefore use one cell from each of these two sets. Since the two sets have sizes three and five, respectively, this is impossible.
\end{proof}

Together, Propositions~\ref{further:through-eight}, \ref{further:nine}, \ref{further:ten}, and~\ref{further:eleven} establish the additional properties stated in Remark~\ref{rem:further}.

\end{document}